\documentclass[12pt]{amsart}

\newtheorem{theorem}{Theorem}[section]
\newtheorem{lemma}[theorem]{Lemma}
\newtheorem{corollary}[theorem]{Corollary}
\newtheorem{proposition}[theorem]{Proposition}
\newtheorem{remark}[theorem]{Remark}
\newtheorem{definition}[theorem]{Definition}
\newtheorem{conjecture}[theorem]{Conjecture}

\newcommand{\ncom}{\newcommand}
\ncom{\rar}{\rightarrow}
\ncom{\lrar}{\longrightarrow}
\ncom{\ov}{\overline}
\ncom{\m}{\mbox}
\ncom{\sta}{\stackrel}
\ncom{\comx}{{\mathbb C}}
\ncom{\Z}{{\mathbb Z}}
\ncom{\Q}{{\mathbb Q}}
\ncom{\R}{{\mathbb R}}
\ncom{\G}{{\mathbb G}}
\ncom{\al}{\alpha}
\ncom{\p}{{\mathbb P}}
\ncom{\E}{{\mathbb E}}
\ncom{\N}{{\mathbb N}}
\ncom{\K}{{\mathbb K}}
\ncom{\Le}{{\mathbb L}}
\ncom{\A}{{\mathbb A}}
\ncom{\B}{{\mathbb B}}
\ncom{\F}{{\mathbb F}}
\ncom{\C}{{\mathbb C}}
\ncom{\f}{\frac}
\ncom{\cA}{{\mathcal A}}
\ncom{\cX}{{\mathcal X}}
\ncom{\cO}{{\mathcal O}}
\ncom{\cW}{{\mathcal W}}
\ncom{\cL}{{\mathcal L}}
\ncom{\cP}{{\mathcal P}}
\ncom{\cH}{{\mathcal H}}
\ncom{\cS}{{\mathcal S}}
\ncom{\cM}{{\mathcal M}}
\ncom{\cC}{{\mathcal C}}
\ncom{\cT}{{\mathcal T}}
\ncom{\cF}{{\mathcal F}}
\ncom{\cN}{{\mathcal N}}
\ncom{\cJ}{{\mathcal J}}
\ncom{\cV}{{\mathcal V}}
\ncom{\cZ}{{\mathcal Z}}
\ncom{\cU}{{\mathcal U}}
\ncom{\cSU}{{\mathcal S \mathcal U}}
\ncom{\cG}{{\mathcal G}}
\ncom{\cQ}{{\mathcal Q}}
\ncom{\cR}{{\mathcal R}}
\ncom{\cE}{{\mathcal E}}

\begin{document}
\baselineskip=16pt

\title[Chow--K\"unneth decompositions]{Chow--K\"unneth decomposition
 for special varieties}

\author[J. N. Iyer]{Jaya NN Iyer}
\author[S. M\"uller--Stach]{Stefan M\"uller--Stach}
\address{The Institute of Mathematical Sciences, CIT
Campus, Taramani, Chennai 600113, India}
\email{jniyer@imsc.res.in}

\address{Mathematisches Institut der Johannes Gutenberg Universit\"at
 Mainz,
Staudingerweg 9, 55099 Mainz, Germany}
\email{mueller-stach@mathematik.uni-mainz.de}

\footnotetext{Mathematics Classification Number: 14C25, 14D05, 14D20, 14D21}

\begin{abstract}In this paper we investigate Murre's conjecture on the Chow--K\"unneth decomposition for two classes of examples. We look at the universal families
of smooth curves over spaces which dominate the moduli space $\cM_g$, in genus at most $8$ and show existence of a Chow--K\"unneth decomposition. The second class of examples include the representation varieties of a finitely generated group with one relation. This is done in the setting of equivariant cohomology and equivariant Chow groups to get equivariant Chow--K\"unneth decompositions. 
\end{abstract}
\maketitle


\setcounter{tocdepth}{1}
\tableofcontents

\section{Introduction}

Suppose $X$ is a nonsingular projective variety defined over the
 complex numbers. We consider the rational Chow group
$CH^i(X)_\Q=CH^i(X)\otimes \Q$ of algebraic cycles of codimension $i$
 on $X$. The conjectures of S. Bloch and A. Beilinson
predict a finite descending filtration $\{F^jCH^i(X)_\Q\}$ on
 $CH^i(X)_\Q$ and satisfying certain compatibility conditions.
A candidate for such a filtration has been proposed by J. Murre and he
 has made the following conjecture \cite{Mu2},

\textbf{Murre's conjecture}: The motive $(X,\Delta)$ of $X$ has a Chow-K\"unneth decomposition:
  $$\Delta= \sum_{i=0}^{2d}\pi_i \in CH^d(X\times X)\otimes \Q$$
such that $\pi_i$ are orthogonal projectors.

These projectors give a candidate for a filtration of the rational Chow groups, see \S \ref{Murrefiltration}.

This conjecture is known to be true for curves, surfaces and
a product of a curve and surface \cite{Mu1}, \cite{Mu3}.
A variety $X$ is known to have a Chow--K\"unneth decomposition if
$X$ is an abelian variety/scheme \cite{Sh},\cite{De-Mu}, a
uniruled threefold \cite{dA-Mul}, universal families over modular varieties
\cite{Go-Mu}, \cite{GHM2} and the universal family over one Picard modular
 surface \cite{MM}, where a partial set of projectors are found. Finite group quotients (maybe singular) of an
abelian variety also satisfies the above conjecture \cite{Ak-Jo}. Furthermore, for some varieties with a nef
 tangent bundle, Murre's conjecture is proved in
 \cite{Iy}.
A criterion for existence of such a decomposition is also given in
 \cite{Sa}. Some other examples are also listed in \cite{Gu-Pe}.

Gordon-Murre-Hanamura \cite{GHM2}, \cite{Go-Mu} obtained  Chow--K\"unneth projectors for universal families  over modular 
varieties. 
 Hence it is natural to ask if the universal families over
the moduli space of curves of higher genus also admit a Chow--K\"unneth decomposition. 
In this paper, we investigate the existence of Chow--K\"unneth decomposition for families of 
smooth curves over spaces which closely approximate the moduli spaces of curves  $\cM_g$ of genus at most $8$, see \S \ref{curves}. The other examples that we
look at are the representation varieties of finitely generated groups with one relation, see \S \ref{repvarieties}. 

In both the above class of examples, we take into account the non-trivial 
action of a linear algebraic group $G$ acting on the spaces. This gives rise to
 the equivariant cohomology and equivariant Chow groups, which were introduced
 and studied by Borel, Totaro, Edidin-Graham \cite{Borel}, \cite{Totaro}, \cite{EdidinGraham}.
Hence it seems natural to formulate Murre's conjecture with respect to the 
cycle class maps between the rational equivariant Chow groups and the rational equivariant cohomology, see \S \ref{equivariantMurre}. Since in concrete examples, good quotients of non-compact varieties exist, it became necessary to extend Murre's conjecture for non-compact smooth varieties, by taking only the bottom weight cohomology $W_iH^i(X,\Q)$ (see \cite{D3}), into consideration. This is weaker than 
the formulation done in \cite{BE}. For our purpose though, it suffices to look at this weaker formulation. We then construct a category of equivariant Chow motives, fixing an algebraic group $G$ (see  \cite{Ak-Jo}, for a category of motives of quotient varieties, under a finite group action).

With this formalism, we show (see \S \ref{familiesofcurves});

\begin{theorem}
The equivariant Chow motive of a universal family of smooth curves 
$\cX\rar U$ over spaces $U$ which dominate the moduli space of curves $\cM_g$, for $g\leq 8$, admits an equivariant Chow--K\"unneth decomposition, for a suitable
linear algebraic group $G$ acting non-trivially on $\cX$.
\end{theorem}

Whenever smooth good quotients exist under the action of $G$, then the equivariant Chow-K\"unneth projectors actually correspond to the absolute Chow--K\"unneth
 projectors for the quotient varieties. In this way, we get orthogonal projectors for universal families over spaces which closely approximate the moduli spaces $\cM_g$, when $g$ is at most $8$.
  
Similarly, the second class of examples we look at are the representation varieties
$R(\Gamma_g,GL_n)$. Here $\Gamma_g$ is a freely generated group, with $2g$ generators
 and with one relation, and the algebraic group $GL_n$ acts by conjugation on this 
variety.

We show
\begin{theorem}
The equivariant Chow motive of the representation variety $R(\Gamma_g,GL_n)$ admits an equivariant Chow--K\"unneth decomposition for the conjugation action by $GL_n$.
\end{theorem}

The other examples that admit a Chow--K\"unneth decomposition are the Fano varieties of $r$-dimensional planes contained in a general complete intersection in a projective space, see Corollary \ref{fano}.
  
The proofs involve classification of curves in genus at most $8$ by Mukai \cite{Mukai},\cite{Mukai2} with respect to embeddings as complete intersections
in homogeneous spaces. This allows us to use Lefschetz theorem and construct orthogonal projectors.

{\Small
\textbf{Acknowledgements}: The first named author thanks the Math Department of Mainz, for their hospitality and financial support during a visit in July 2006 and Sept 2007, when this work was carried out.
}


\section{Preliminaries}

The category of nonsingular projective varieties over $\comx$ will be
denoted by $\cV$. Let $CH^i(X)_\Q=CH^i(X)\otimes \Q$ denote the rational Chow group of
codimension $i$ algebraic cycles modulo rational equivalence.

Suppose $X,Y\in Ob(\cV)$ and $X=\cup X_i$ be a decomposition into
connected components $X_i$ and $d_i=\m{dim }X_i$.
Then $\m{Corr}^r(X,Y)= \oplus_i CH^{d_i+r}(X_i\times Y)_\Q$ is called a
space of correspondences of degree $r$ from $X$ to
$Y$.

Suppose $X$ is a nonsingular projective variety over $\comx$ of
dimension $d$. Let $\Delta\subset X\times X$ be the diagonal.
Consider the K\"unneth decomposition of $\Delta$ in the Betti Cohomology:
$$\Delta = \oplus_{i=0}^{2d}\pi_i^{hom}$$
where $\pi_i^{hom}\in H^{2d-i}(X)\otimes H^i(X)$.

\begin{definition}
The motive of $X$ is said to have K\"unneth decomposition if each of
the classes $\pi_i^{hom}$ are algebraic i.e., $\pi_i^{hom}$ is the image of an algebraic cycle $\pi_i$
under the cycle class map from the rational Chow groups to the Betti cohomology.

\end{definition}

\begin{definition}
The motive of $X$ is said to have a Chow--K\"unneth decomposition if
 each of the
classes $\pi_i^{hom}$ is algebraic and are orthogonal projectors,
 i.e.,
$\pi_i\circ \pi_j=\delta_{i,j}\pi_i$.
\end{definition}

\begin{lemma}\label{le.-1}
If $X$ and $Y$ have a Chow--K\"unneth decomposition then $X\times Y$
 also
has a Chow--K\"unneth decomposition.
\end{lemma}
\begin{proof}
If $\pi_i^X$ and $\pi_j^Y$ are the Chow--K\"unneth components for
 $h(X)$
and $h(Y)$ respectively then
$$
\pi_i^{X\times Y}=\sum_{p+q=i}\pi_p^X\times \pi_q^Y
 \,\in\,CH^*(X\times Y \times X\times Y)_\Q
$$
are the Chow--K\"unneth components for $X\times Y$. Here the product
$\pi_p^X\times \pi_q^Y$ is taken after identifying $X\times
Y\times X\times Y \simeq X\times X \times Y\times Y$.
\end{proof}

\subsection{Murre's conjectures}\label{Murrefiltration}

J. Murre \cite{Mu2}, \cite{Mu3} has made the following conjectures for any smooth projective variety $X$. 

(A) The motive $h(X):=(X,\Delta_X)$ of $X$ has a Chow-K\"unneth decomposition:
  $$\Delta_X= \sum_{i=0}^{2n}\pi_i \in CH^n(X\times X)\otimes \Q$$
such that $\pi_i$ are orthogonal projectors.
 
(B) The correspondences $\pi_0,\pi_1,...,\pi_{j-1},\pi_{2j+1},...,\pi_{2n}$ act as zero on $CH^j(X)\otimes \Q$.

(C) Suppose $$F^rCH^j(X)\otimes \Q= Ker \pi_{2j}\cap Ker \pi_{2j-1}\cap...\cap Ker \pi_{2j-r+1}.$$
 Then the filtration $F^\bullet$ of $CH^j(X)\otimes \Q$ is independent of the choice of the projectors $\pi_i$.

(D) Further, $F^1CH^i(X)\otimes \Q= (CH^i(X)\otimes \Q)_{hom}$, the cycles which are homologous to zero.
\vskip .3cm

In \S \ref{Murreequiv}, we will extend (A) in the setting of equivariant
Chow groups.


\section{Equivariant Chow groups and equivariant Chow motives}


In this section, we recall some preliminary facts on the equivariant groups to formulate Murre's conjectures for a smooth variety $X$ of dimension $d$, which is equipped with an action by a linear reductive algebraic group $G$.
The equivariant groups and their properties that we recall below were defined by Borel, Totaro, Edidin-Graham, Fulton \cite{Borel},\cite{Totaro},\cite{EdidinGraham}, \cite{Fulton}. 

\subsection{Equivariant cohomology $H^i_G(X,\Z)$ of $X$}

Suppose $X$ is a variety with an action on the left by an algebraic group $G$. Borel defined the equivariant cohomology $H^*_G(X)$ as follows.
There is a contractible space $EG$ on which $G$ acts freely (on the right) with quotient $BG:=EG/G$. Then form the space
$$
EG\times_G X:=EG\times X/(e.g,x)\simeq (e,g.x).
$$
In other words, $EG\times_G X$ represents the (topological) quotient stack $[X/G]$.

\begin{definition} The equivariant cohomology of $X$ with respect to $G$ is the ordinary singular cohomology of  $EG\times_G X$:
$$
H^i_G(X)= H^i(EG\times_G X).
$$
\end{definition}
For the special case when $X$ is a point, we have
$$
H^i_G(point)= H^i(BG)
$$
For any $X$, the map $X\rar point$ induces a pullback map $ H^i(BG)\rar H^i_G(X)$. Hence the equivariant cohomology of $X$ has the structure of a $ H^i(BG)$-algebra, at least when $ H^i(BG)=0$ for odd $i$.

\subsection{Equivariant Chow groups $CH^i_G(X)$ of $X$}\cite{EdidinGraham}\label{equivchow}
As in the previous subsection, let $X$ be a smooth variety, equipped with a left $G$ action. Here $G$ is an affine algebraic group of dimension $g$. Choose an $l$-dimensional representation $V$ of $G$ such that
$V$ has an open subset $U$ on which $G$ acts freely and whose complement has codimension more than $n-i$. The diagonal action on $X\times U$ is also free, so there is a quotient in the category of algebraic spaces.
Denote this quotient by $X_G:=(X\times U)/G$.

\begin{definition} The $i$-th equivariant Chow group $CH^G_i(X)$ is the usual Chow group $CH_{i+l-g}(X_G)$.
\end{definition}

\begin{proposition}
The equivariant Chow group $CH^G_i(X)$ is independent of the representation $V$, as long as $V-U$ has sufficiently high codimension.
\end{proposition}
\begin{proof}
See \cite[Definition-Proposition 1]{EdidinGraham}.
\end{proof}

If $Y\subset X$ is an $m$-dimensional subvariety which is invariant under the $G$-action, then it has a $G$-equivariant fundamental class
$[Y]_G \in CH^G_m[X]$. More generally if $V$ is an $l$-dimensional representation of $G$ and $S\subset X\times V$ is an $m+l$-dimensional subvariety which is invariant under $G$-action, then $S$ has a $G$-equivariant fundamental class $[S]_G\in CH^G_m(X)$.

\begin{proposition}
If $\alpha\in CH^G_m(X)$ then there exists a representation $V$ such that $\alpha= \sum a_i[S_i]_G$, for some $G$-invariant subvarieties $S_i$
of $X\times V$.
\end{proposition} 
\begin{proof}
See \cite[Proposition 1]{EdidinGraham}.
\end{proof}

\subsection{Functoriality properties}

Suppose $f:X\rar Y$ is a $G$-equivariant morphism.
Let $\cS$ be one of the following properties of schemes or algebraic spaces: proper, flat, smooth, regular embedding or l.c.i.

\begin{proposition}
If $f:X\rar Y$ has property $\cS$, then the induced map
$f_G:X_G\rar Y_G$ also has property $\cS$.
\end{proposition}
\begin{proof}
See \cite[Proposition 2]{EdidinGraham}.
\end{proof}

\begin{proposition}
Equivariant Chow groups have the same functoriality as ordinary Chow groups for equivariant morphisms with property $\cS$.
\end{proposition}
\begin{proof}
See \cite[Proposition 3]{EdidinGraham}.
\end{proof}

If $X$ and $Y$ have $G$-actions then there are exterior products
$$
CH^G_i(X)\otimes CH^G_j(Y)\rar  CH^G_{i+j}(X\times Y).
$$

In particular, if $X$ is smooth then there is an intersection product
on the equivariant Chow groups which makes $\oplus_j CH^G_j(X)$ 
into a graded ring.

\subsection{Cycle class maps} \cite[\S 2.8]{EdidinGraham}

Suppose $X$ is a complex algebraic variety and $G$ is a complex algebraic group. The equivariant Borel-Moore homology $H^G_{BM,i}(X)$ is the Borel-Moore homology  $H_{BM,i}(X_G)$, for $X_G=X\times_G U$. This is independent of the representation as long as $V-U$ has sufficiently large codimension. This gives a cycle class map,
$$
cl^i:CH^i_G(X)\rar H_{BM,2i}(X,\Z).
$$
compatible with usual operations on equivariant Chow groups.
Suppose $X$ is smooth of dimension $d$ then $X_G$ is also smooth. In this case the Borel-Moore cohomology $H_{BM,2i}(X,\Z)$ is dual to 
$H^{2d-i}(X_G)= H^{2d-i}(X\times_G U)$.

This gives the cycle class maps
\begin{equation}\label{cyclemap}
cl^i:CH^i_G(X)\rar H^{2i}_G(X,\Z).
\end{equation}

There are also maps from the equivariant groups to the usual groups:
\begin{equation}\label{liftcoh}
H^i_G(X,\Z)\rar H^i(X,\Z)
\end{equation}
and
\begin{equation}\label{liftchow}
CH^i_G(X)\rar CH^i(X).
\end{equation}

\subsection{Weight filtration $W_.$ on $H^i_G(X,\Z)$} 
In this paper, we assign only the bottom weight $W_{i}$ of the equivariant cohomology in the simplest situation. Consider a smooth variety $X$ equipped with a left $G$ action as above. Suppose there is a smooth compactification $\ov{X}$ of $X$
such that the action of $G$ extends to an action on $\ov{X}$.
In this situation, there is a localization sequence, given by the inclusion $j:X\hookrightarrow \ov{X}$:
$$
\rar H^i_G(\ov{X})\sta{j^*}{\rar} H^i_G({X})\rar.
$$
The image under the map $j^*$ is defined to the bottom weight $W_iH^i_G({X})$ cohomology \cite{D3}.
In particular, if $X$ is itself a complete smooth variety then
 $W_iH^i_G({X}) = H^i_G({X})$.

\subsection{Equivariant Chow motives and the category of equivariant Chow motives}
When $G$ is a finite group then a category of Chow motives for (maybe singular) quotients of varieties under the $G$-action was constructed in \cite{Ak-Jo}. More generally, we consider the following situation, taking into account the equivariant cohomology and the equivariant rational Chow groups, which does not seem to have been considered before.
 
Fix an affine complex algebraic group $G$.
Let $\cV_G$ be the category whose objects are complex smooth varieties with a $G$-action and the morphisms are $G$-equivariant morphisms.

For any $X,Y,Z\in Ob(\cV_G)$, consider the projections

$$
X\times Y\times Z \sta{p_{XY}}{\lrar} X\times Y,
$$
$$
X\times Y\times Z \sta{p_{YZ}}{\lrar} Y\times Z,
$$
$$
X\times Y\times Z \sta{p_{XZ}}{\lrar} X\times Z.
$$
which are $G$-equivariant.

Let $d$ be the dimension of $X$.
The group of correspondences from $X$ to $Y$ of degree $r$ is defined as
$$
\m{Corr}^r_G(X\times Y):=CH^{r+d}_G(X\times Y).
$$
Every $G$-equivariant morphism $X\rar Y$ defines an element in $\m{Corr}^0_G(X\times Y)$, by taking the graph cycle.

For any $f\in \m{Corr}^r_G(X,Y)$ and $g\in \m{Corr}^e_G(Y,Z)$ define their composition
$$
g\circ f\in \m{Corr}^{r+e}_G(X,Z)
$$
by the prescription
$$
g\circ f= p_{XZ*}(p_{XY}^*(f).p_{YZ}^*(g)).
$$
 
This gives a linear action of correspondences on the equivariant Chow groups
$$
\m{Corr}^r_G(X,Y)\times CH^s_G(X)_\Q \lrar CH^{r+s}_G(Y)_Q
$$
$$ 
(\gamma,\al) \mapsto p_{Y*}(p_X^*\al.\gamma)
$$
for the projections $p_X:X\times Y\lrar X,\, p_Y:X\times Y\lrar Y$.

The category of pure equivariant $G$-motives with rational coefficients is denoted by $\cM^+_G$. The objects of $\cM^+_G$ are triples
$(X,p,m)_G$, for $X\in Ob(\cV_G),\,p\in \m{Corr}^0_G(X,X)$ is a projector, i.e., $p\circ p=p$ and $m\in \Z$.
The morphisms between the objects $(X,p,m)_G,(Y,q,n)_G$ in $\cM^+_G$ are given by the correspondences
$f\in \m{Corr}^{n-m}_G(X,Y)$ such that $f\circ p=q\circ f=f$. The composition of the morphisms is the compositions of the
correspondences.
This category is pseudoabelian and $\Q$-linear. Furthermore, it is a tensor category defined by 
$$
(X,p,m)_G\otimes (Y,q,n)_G=(X\times Y,p\otimes q,m+n)_G.
$$ 

The object $(\m{Spec}\, \comx,id,0)_G$ is the unit object and the Lefschetz motive $\Le$ is the object $(\m{Spec}\, \comx,id,-1)_G$. Here $\m{Spec}\, \comx$ is taken with a trivial $G$-action.
The Tate twist of a $G$-motive $M$ is $M(r):= M\otimes \Le^{\otimes -r}=(X,p,m+r)_G$.

\begin{definition}
The theory of equivariant Chow motives (\cite{Sc}) provides a functor
$$
h:\cV_G\lrar \cM^+_G.
$$
For each $X\in Ob(\cV_G)$ the object $h(X)=(X,\Delta,0)_G$ is called the equivariant Chow motive of $X$. Here $\Delta$ is the class of the diagonal in $CH^*(X\times X)_\Q$, which is $G$-invariant for the diagonal action on $X\times X$ and hence lies in $\m{Corr}^0_G(X,X)=CH^*_G(X\times X)_\Q$.
\end{definition}

\section{Murre's conjectures for the equivariant Chow motives}\label{Murreequiv}


Suppose $X$ is a complex smooth variety of dimension $d$, equipped with a $G$-action.
Consider the product variety $X\times X$ together with the diagonal action
of the group $G$.

The cycle class map
\begin{equation}\label{cycleclassmap}
cl^d:CH^d(X\times X)_\Q \rar H^{2d}(X\times X,\Q).
\end{equation}
actually maps to the weight $2d$ piece $W_{2d}H^{2d}(X\times X,\Q)$
of the ordinary cohomology group. 

Applying this to the spaces $X\times U$, for open subset $U\subset V$ as in \S \ref{equivchow}, \eqref{cycleclassmap} holds for the equivariant groups as well and there are cycle class maps:

\begin{equation}\label{cycleclassmapG}
cl^d:CH^d_G(X\times X)_\Q \rar W_{2d}H^{2d}_G(X\times X,\Q).
\end{equation}

\begin{lemma}
The image of the diagonal cycle $[\Delta_X]$ under the cycle class map $cl^d$
lies in the subspace 
$$
\bigoplus_{i}W_{2d-i}H^{2d-i}_G(X)\otimes W_iH^i_G(X)
$$
of $W_{2d}H^{2d}_G(X\times X,\Q)$.
\end{lemma}
\begin{proof}
If $X$ is a compact smooth variety then we notice that the weight $2d$ piece 
coincides with the cohomology group $H^{2d}(X\times X,\Q)$ and by the 
K\"unneth formula for products the statement follows in the usual cohomology.
Suppose $X$ not compact. Since the diagonal cycle is $G$-invariant, by \eqref{cycleclassmap}, the image of $[\Delta_X]$ lies in $W_{2d}H^{2d}(X\times X,\Q)$.  Choose a smooth compactification $\ov{X}$ of $X$ and consider the commutative diagram:
\begin{eqnarray*}
\bigoplus_{i}H^{2d-i}(\ov{X})\otimes H^i(\ov{X}) &\sta{\simeq}{\rar} & H^{2d}(\ov{X},\Q)\\
 \downarrow\,\,\,\,\,\,\,\,\, & &\downarrow \\
\bigoplus_{i}W_{2d-i}H^{2d-i}(X)\otimes W_iH^i(X) &\sta{\simeq}{\rar} &W_{2d}H^{2d}(X,\Q)
\end{eqnarray*}
The vertical arrows are surjective maps, defined by the localization.
The assertion now follows from the above diagram in the bottom weight of the ordinary cohomology group of any smooth variety.
In particular, it is true in the bottom weights of the ordinary cohomology groups of the smooth variety $X\times U$, for any open subset $U\subset V$ of large complementary codimension and $V$ is a $G$-representation. 
 But this is essentially the bottom weights of the equivariant cohomology group of $X$.

\end{proof}

Denote the decomposition of the $G$-invariant diagonal cycle
\begin{equation}\label{equiprojectors}
\Delta_X= \oplus_{i=0}^{2d}\pi_i^G \,\in \,W_{2d}H^{2d}_G(X,\Q)
\end{equation}
such that $\pi_i^G$ lies in the space $W_{2d-i}H^{2d-i}_G(X)\otimes W_iH^i_G(X)$.

\begin{definition}
Suppose $X$ is a smooth variety equipped with a $G$-action.
The equivariant Chow motive $(X,\Delta_X)_G$ of $X$ is said to have 
an $\textbf{equivariant K\"unneth decomposition}$ if the classes $\pi_i^G$ are algebraic, i.e., they have a lift in the equivariant Chow group $CH^d_G(X\times X)_\Q$. Furthermore, if the lifts are orthogonal projectors then we say that the equivariant Chow motive of $X$ has an \textbf{equivariant Chow--K\"unneth decomposition}.
\end{definition}

We can now extend Murre's conjecture to smooth varieties with a $G$-action, as follows.

\begin{conjecture}\label{equivariantMurre}
Suppose $X$ is a smooth variety with a $G$-action. Then the equivariant
Chow motive $(X,\Delta_X)_G$ of $X$ has an equivariant Chow--K\"unneth decomposition.
\end{conjecture}  

In particular, if the action of $G$ is trivial then we can extend Murre's conjecture to a (not necessarily compact) smooth variety, by taking only the bottom weight cohomology $W_iH^i(X)$ of the ordinary cohomology. This is weaker than obtaining projectors for the ordinary cohomology.
We remark a projector $\pi_1$ in the case of quasi--projective varieties has been constructed by Bloch and Esnault~\cite{BE}.

\section{Families of curves}\label{curves}


Our goal in this paper is to find (explicit) absolute Chow--K\"unneth 
decomposition for the universal families of curves over close approximations of the moduli space of smooth curves of small genus.
We begin with the case of families of plane curves.  

\subsection{Family of plane curves}\label{family}

There is an easy method to obtain a family of plane curves by blowing up a line on a hypersurface, which is illustrated below. We show that this family has an absolute Chow--K\"unneth
decomposition. This motivates the other examples we look in this paper.
 
Suppose $X\subset \p^n$ is a nonsingular hypersurface of degree $d$ and
 containing a line.
If $n\geq 3$ and $d\leq n+1$ then $X$ always contains a line
 \cite{Kol}. Let 
$L$ be any line on $X$ and we assume that $L$ is not contained in any
 plane $P\subset X$.

Projecting from $L$, we obtain a diagram (A),

\begin{eqnarray*}
X'& \sta{f}{\lrar} & X \\
\downarrow \!p & & \\
\p^{n-2}
\end{eqnarray*}
such that $X'$ is the blow-up of $X$ along $L$ and $p:X'\lrar \p^{n-2}$
 is a family of curves of degree $d-1$, with the generic fiber smooth.
This family can be described as follows.

Consider the morphism 
$$
\pi':\p^{n-2}\lrar G(2,\p^n),\,x\mapsto P=<x,L>
$$
which is clearly an embedding.
 Here $G(2,\p^n)$ denotes the Grassmanian of $2$--planes contained in
 $\p^n$ and $P=<x,L>$ is the plane spanned by
$x$ and the line $L$. 

Consider the morphism
\begin{equation}\label{map}
\eta:\m{Image }(\pi')\lrar |\cO_{\p^2}(d-1)|,\, P\mapsto C
\end{equation}
where $P\cap X=L+C$ where $C\subset P$ is a curve of degree $d-1$.
Then $X'\rar \p^{n-2}$ is the pullback of the universal plane curve over 
$|\cO_{\p^2}(d-1)|$.

\begin{lemma}\label{absCK}
The family 
$$
p: X'\lrar \p^{n-2}.
$$
of projective plane curves has an absolute Chow--K\"unneth
 decomposition.
\end{lemma}
\begin{proof}
We first notice that a smooth hypersurface $X\subset \p^n$ of degree
 $d$ has a Chow--K\"unneth decomposition. Indeed, the cohomology of $X$ is
 algebraic except in the middle dimension $H^{n-1}(X,\Q)$. By the
 Lefschetz Hyperplane section theorem, the algebraic cohomology
 $H^{2j}(X,Q)$, $j\neq n-1$, is generated by the hyperplane section $H^j$.
So the projectors are simply
$$
\pi_r:= \f{1}{d}.H^{n-1-r}\times H^r \,\in\, CH^{n-1}(X\times X)_\Q
$$
for $r\neq n-1$.
We can now take $\pi_{n-1}:= \Delta_X- \sum_{r, r\neq n-1} \pi_r$. 
This gives a complete set of orthogonal projectors and a
 Chow--K\"unneth decomposition for $X$.
Since $X'\rar X$ is a blow-up along a line, the new cohomology is again
 algebraic, by the blow-up formula. Similarly we get a Chow--K\"unneth
 decomposition for $X'$ (see also \cite[Lemma 2]{dA-Mu2} for blow-ups).
\end{proof}

The above lemma can be generalized to the following situation.

\begin{lemma}\label{simpleprojectors}
Suppose $Y$ is a smooth projective variety of dimension $r$ over $\comx$ which
has only algebraic cohomology groups $H^i(Y)$ for all $0 \le i \le m$ for some $m < r$. 
Then we can construct orthogonal projectors
$$
\pi_0,\pi_1,...,\pi_m,\pi_{2r-m},\pi_{2r-m+1},...,\pi_{2r}
$$
in the usual Chow group $CH^r(Y\times Y)_\Q$, and where $\pi_{2i}$ acts as $\delta_{i,p}$ on $H^{2p}(Y)$ and $\pi_{2i-1}=0$.
Moreover, if there is an affine complex algebraic group $G$ acting on $Y$, then we can lift the above projectors
in the equivariant Chow group $CH^r_G(Y\times Y)_\Q$ as orthogonal projectors.
\end{lemma}
\begin{proof} See also \cite{dA-Mul, dA-Mu2}. Let $H^{2p}(Y)$ be generated by cohomology classes of cycles $C_1,\ldots,C_s$ and $H^{2r-2p}(Y)$  
be generated by cohomology classes of cycles $D_1,\ldots,D_s$. We denote by $M$ the intersection matrix with entries
$$
M_{ij}= C_i \cdot D_j \in \Z.
$$ 
After base change and passing to $\Q$--coefficients we may assume that $M$ is diagonal, 
since the cup--product $H^{2p}(Y,\Q) \otimes H^{2r-2p}(Y,\Q) \to \Q$ is non--degenerate. We define the projector $\pi_{2p}$ as 
$$
\pi_{2p}=\sum_{k=1}^s \frac{1}{M_{kk}} D_k \times C_k. 
$$
It is easy to check that $\pi_{2p\,*}(C_k)=D_k$. Define $\pi_{2r-2p}$ as the adjoint, i.e., transpose of $\pi_{2p}$.
Via the Gram--Schmidt process from linear algebra we can successively make all projectors orthogonal.
\end{proof}
 
Suppose $X\subset \p^n$ is a smooth complete intersection of multidegree $d_1\leq d_2\leq ...\leq d_s$. Let $F_r(X)$ be the variety of $r$-dimensional planes
contained in $X$. Let $\delta:=\min\{(r+1)(n-r)- {d+r \choose r},n-2r-s\}$.

\begin{corollary}\label{fano}
 If $X$ is general then $F_r(X)$ is a smooth projective variety of dimension $\delta$ and it has an absolute Chow--K\"unneth decomposition.
\end{corollary}
\begin{proof}
The first assertion on the smoothness of the variety $F_r(X)$ is well--known, see \cite{Altman}, \cite{ELV}, \cite{De-Ma}.
For the second assertion, notice that $F_r(X)$ is a subvariety of the Grassmanian $G(r, \p^n)$ and is the zero set of a section of a vector bundle. Indeed, let $S$ be the tautological bundle on $G(r,\p^n)$. Then a section of $\oplus_{i=1}^s Sym^{d_i}H^0(\p^n,\cO(1))$ induces a section of the vector bundle $\oplus_{i=1}^s Sym^{d_i}S^*$ on $G(r,\p^n)$.
Thus, $F_r(X)$ is the zero locus of the section of the $\bigoplus_{i=1}^s Sym^{d_i}S^*$ induced by the equations defining the complete intersection $X$.
A Lefschetz theorem is proved in \cite[Theorem 3.4]{De-Ma}:
$$
H^i(G(r,\p^n),\Q)\rar  H^i(F_r(X),\Q)
$$
is bijective, for $i\leq \delta-1$.
We can apply Lemma \ref{simpleprojectors} to get the orthogonal projectors
in all degrees except in the middle dimension. The projector corresponding to the middle dimension can be gotten by subtracting the sum of these projectors from the diagonal class.

\end{proof}

\begin{corollary}
Suppose $X\subset \p^n$ is a smooth projective variety. Let $Y\subset\p^n$ be any other projective variety such that $X\cap Y$ is smooth and irreducible. If $X$ has only algebraic cohomology, except perhaps in the middle dimension, then $X\cap Y$ has an absolute Chow--K\"unneth decomposition.
\end{corollary}
\begin{proof}
Since the class of $Y$ in the cohomology of $\p^n$ is a linear section, the Lefschetz theorem holds for the restriction map
$$
H^i(X,\Q)\sta{\cap Y}{\rar} H^i(X\cap Y,\Q).
$$ 
Now we can apply Lemma \ref{simpleprojectors} to deduce our claim.
\end{proof}

\begin{corollary}
Suppose $X\subset \p^n$ is a smooth projective variety of dimension $d$. Let
$r=2d-n$.
 Then we can construct orthogonal projectors
$$
\pi_0,\pi_1,...,\pi_r,\pi_{2d-r},\pi_{2d-r+1},...,\pi_{2d}.
$$
\end{corollary}
\begin{proof}
Barth \cite{Ba} has proved a Lefschetz theorem for higher codimensional subvarieties in projective spaces:
$$
H^i(\p^n,\Q)\rar  H^i(X,\Q)
$$
is bijective if $i\leq 2d-n$ and is injective if $i=2d-n+1$.
The claim now follows from Lemma \ref{simpleprojectors}.
\end{proof}

\begin{remark}
The above corollary says that if we can embed a variety $X$ in a low dimensional projective space then we get at least a partial set 
of orthogonal projectors. A conjecture of Hartshorne's says that any codimension two subvariety of $\p^n$ for $n \ge 6$ is a complete intersection. 
This gives more examples for subvarieties with several algebraic cohomology groups.
\end{remark}

\subsection{Chow--K\"unneth decomposition for the universal plane curve}

We want to find explicit equivariant Chow--K\"unneth projectors for the universal plane curve of degree $d$. Let $d\geq 1$ and consider the linear system
 $\p=|\cO_{\p^2}(d)|$ and the universal plane curve

\begin{eqnarray*}
\cC & \subset & \p^2 \times \p \\
\downarrow && \\
\p. &&
\end{eqnarray*}

Furthermore, we notice that the general linear group $G:=GL_3(\comx)$ acts on 
$\p^2$ and hence acts on the projective space $\p=|\cO_{\p^2}(d)|$. This 
gives an action on the product space $\p^2 \times \p$ and leaves the universal plane curve $\cC\subset  \p^2 \times \p$ invariant under the $G$-action.

\begin{lemma}
The variety $\cC$ has an absolute Chow--K\"unneth decomposition and an absolute equivariant Chow--K\"unneth decomposition.
\end{lemma}
\begin{proof}
We observe that $\cC\subset  \p^2 \times \p$ is a smooth hypersurface
 of degree $d$ with variables in $\p^2$ whose coefficients are polynomial
 functions on $\p$. Notice that $\p^2\times \p$ has a Chow--K\"unneth
 decomposition and Lefschetz
theorems hold for the embedding $\cC\subset \p^2 \times \p$. Now we can
 repeat the arguments from Lemma \ref{simpleprojectors} to get an absolute
 Chow--K\"unneth decomposition and absolute equivariant Chow--K\"unneth decomposition, for the variety $\cC$.
\end{proof}

\subsection{Families of curves contained in homogeneous spaces}\label{familiesofcurves}
We notice that when $d=3$ in the previous subsection, the family of plane 
cubics restricted to the loci of stable curves is a complete family of genus one stable curves. If $d\geq 4$, then the above family of plane curves is no longer a complete family of genus $g$ curves. Hence
to find families which closely approximate over the moduli spaces of stable curves, we need to look for curves embedded as complete intersections in other simpler looking varieties.
For this purpose, we look at the curves embedded in special Fano varieties of small genus $g\leq 8$, which was studied by S. Mukai \cite{Mukai}, \cite{Mukai2}, \cite{Mukai3}, \cite{Mukai5} and Ide-Mukai \cite{Mukai4}.

We recall the main result that we need.

\begin{theorem}\label{mukai}
Suppose $C$ is a generic curve of genus $g\leq 8 $. Then $C$ is a complete intersection in a smooth projective variety which has only algebraic cohomology.
\end{theorem}
\begin{proof}
This is proved in \cite{Mukai}, \cite{Mukai2}, \cite{Mukai3}, \cite{Mukai4} and \cite{Mukai5}.
The below classification is for the generic curve.

When $g\leq 5$ then it is well-known that the generic curve is a linear section of a Grassmanian.

When $g=6$ then a curve has finitely many $g^1_4$ if and only if it is a complete intersection of a Grassmanian and a smooth quadric , see \cite[Theorem 5.2]{Mukai3}.

When $g=7$ then a curve is a linear section of a $10$-dimensional spinor variety $X\subset \p^{15}$ if and only if it is non-tetragonal, see \cite[Main theorem]{Mukai3}. 

When $g=8$ then it is classically known that the generic curve is a linear section of the grassmanian $G(2,6)$ in its Plucker embedding.

\end{proof}

Suppose $\p(g)$ is the parameter space of linear sections of a Grassmanian or of a spinor variety, which depends on the genus, 
as in the proof of above Theorem \ref{mukai}. $\p(g)$ is a product of projective spaces on which an algebraic group $G$ 
(copies of $PGL_N$) acts. Generic curves are isomorphic, if they are in the same orbit of $G$.

\begin{proposition}\label{generic}
Suppose $\p(g)$ is as above, for $g\leq 8$. Then there is a universal curve
$$
\cC_g\rar \p(g)
$$
such that the classifying (rational) map $\p(g)\rar \cM_g$ is dominant.
The smooth variety $\cC_g$ has an absolute Chow--K\"unneth decomposition and an
absolute equivariant Chow-K\"unneth decomposition for the natural $G$--action mentioned above. 
\end{proposition}
\begin{proof}
The first assertion follows from Theorem \ref{mukai}. For the second assertion 
notice that the universal curve, when $g\leq 8$, is a complete intersection
in $\p(g)\times V$ where $V$ is either a Grassmanian or a spinor variety, which
are homogeneous varieties. In other words, $\cC_g$ is a complete intersection in a space which has only algebraic cohomology. Hence, by Lemma \ref{simpleprojectors}, $\cC_g$ has an absolute Chow--K\"unneth decomposition.
Now a homogeneous variety looks like $V=G/P$ where $G$ is an (linear) algebraic group and
$P$ is a parabolic subgroup. Hence the group $G$ acts on the variety $V$. This
induces an action on the linear system $\p(g)$ and hence $G$ acts on the ambient variety $\p(g)\times V$ and leaves the universal curve $\cC_g$ invariant.
Hence we can again apply Lemma \ref{simpleprojectors} to obtain absolute equivariant Chow--K\"unneth decomposition for $\cC_g$.   
 
\end{proof}

Consider the universal family of curves $\cC_g\rar \p(g)$ as obtained above, 
which are equipped with an action of a linear algebraic group $G$.

Suppose there is an open subset $U_g\subset \p(g)$, with the universal family 
$\cC_{U_g}\rar U_g$, on which $G$ acts freely to form
a good quotient family 
$$ 
Y_g:=\cC_{U_g}/G\rar S_g:=U_g/G.
$$ 
Notice that the classifying map $S_g\rar \cM_g$ is dominant.

\begin{corollary}
The smooth variety $Y_g$ has an absolute Chow--K\"unneth decomposition.
\end{corollary}
\begin{proof}
Consider the localization sequence, for the embedding $j:\cC_{U_g}\times \cC_{U_g}\hookrightarrow \cC_g\times \cC_g$,
$$
CH^d_G(\cC_g\times \cC_g)_\Q \sta{j^*}{\rar} CH^d_G(\cC_{U_g}\times \cC_{U_g})_\Q\rar 0.
$$
Here $d$ is the dimension of $\cC_g$.
Then the map $j^*$ is an equivariant ring homomorphism and transforms orthogonal projectors
to orthogonal projectors.
Similarly there is a commuting diagram between the equivariant cohomologies:
\begin{eqnarray*}
\bigoplus_{i}H^{2d-i}_G(\cC_g)\otimes H^i_G(\cC_g) &\sta{\simeq}{\rar} & H^{2d}_G(\cC_g,\Q)\\
 \downarrow\,\,\,\,\,\,\,\,\, & &\downarrow \\
\bigoplus_{i}W_{2d-i}H^{2d-i}_G(\cC_{U_g})\otimes W_iH^i_G(\cC_{U_g}) & \sta{\simeq}{\rar}&W_{2d}H^{2d}_G(\cC_{U_g},\Q)
\end{eqnarray*}
The vertical arrows are surjective maps mapping onto the bottom weights of the equivariant cohomology groups.
By Proposition \ref{generic}, the variety $\cC_g$ has an absolute equivariant Chow--K\"unneth decomposition.
Hence the images of the equivariant Chow--K\"unneth projectors for the complete smooth variety $\cC_g$, under the morphism $j^*$ give equivariant Chow--K\"unneth projectors for the smooth variety $\cC_{U_g}$.
 
Using \cite{EdidinGraham}, we have the identification of the rational Chow groups
$$
CH^*(Y_g)_\Q\,=\,CH^*_G(\cC_{U_g})_\Q
$$
and
$$
CH^*(Y_g\times Y_g)_\Q\,=\,CH^*_G(\cC_{U_g}\times \cC_{U_g})_\Q.
$$
Furthermore, this respects the ring structure on the above rational Chow groups.
A similar identification also holds for the rational cohomology groups.
This means that the equivariant Chow--K\"unneth projectors for the variety 
$\cC_{U_g}$ correspond to a complete set of absolute Chow--K\"unneth projectors
for the quotient variety $Y_g$.
\end{proof}

\begin{remark}
Since Mukai has a similar classification for the non-generic curves in genus
$\leq 8$, one can obtain absolute equivariant Chow--K\"unneth decomposition
for these special families of smooth curves, by applying the proof of Proposition \ref{generic}. There is also a classification for $K3$-surfaces and in many cases the generic $K3$-surface is obtained as a linear section of a Grassmanian
\cite{Mukai}. Hence we can apply the above results to families of $K3$-surfaces over spaces which dominate the moduli space of $K3$-surfaces.
\end{remark}

\section{Equivariant Chow--K\"unneth decomposition for representation varieties}\label{repvarieties}


In this section, we look at the representation varieties $R(\Gamma,GL_n)$, where
$\Gamma$ is a finitely generated group with one relation. It is significant to consider such varieties because these varieties allow a conjugation action by 
$GL_n$. When $\Gamma$ is the fundamental group of a compact Riemann surface, the good quotients of the semi-stable loci give rise to moduli spaces. Instead of $GL_n$, if we consider the unitary group $U(n)$, then the study of such quotients was initiated by Mumford \cite{Md}, \cite{Md2} and Narasimhan-Seshadri \cite{Na-Se}. 
More generally, when $\Gamma$ is the fundamental group of a complex smooth projective variety then Simpson \cite{Si}, \cite{Si2}, \cite{Si3} studied such quotients and proved the equivalences between the Betti moduli space, the de Rham moduli space and the Dolbeault moduli space. This equivalence is although not algebraic.
Hence the equivariant Chow--K\"unneth projectors for $R(\Gamma,GL_n)$ will correspond to Chow-K\"unneth projectors for the moduli space of finite dimensional irreducible representations, which is a smooth variety.

\subsection{Representation varieties}\cite{RBC}
Let $\Gamma$ be a finitely generated group. For any algebraic group $G$, the set $R(\Gamma, G)$ of all representations $\rho: \Gamma\rar G$ is known to have a natural structure of an algebraic variety. With this structure, this set is called the representation variety of $\Gamma$ in $G$. When $G=GL_n$, this variety is called the variety of $n$-dimensional representations of $\Gamma$ and usually denoted by $R_n(\Gamma)$.

Suppose $\Gamma$ has $2g$ generators and admits at most one relation $r$.
Such groups arise as the fundamental groups of compact Riemann surfaces and of 
abelian varieties or compact complex tori, when there is no relation.
When there are $2g$ generators and one relation, denote the representation variety by $R_n(\Gamma_g)$. When there are $2g$ generators and no relation, denote the representation variety by $C(2g, n):= R(\Z_{2g},GL_n)$, which is usually called the \textit{Commutator variety}. This can also be described as the set of $2g$-tuples of pairwise commuting $2g\times 2g$-matrices;
\begin{equation}\label{character}
C(2g,GL_n):= \{(x_1,x_2,...,x_{2g})\in (GL_n)^{2g}|\,\,x_i.x_j=x_j.x_i, \m{ for all } i,j\leq 2g\}.
\end{equation}
 
Similarly, we can write
\begin{equation}\label{relation}
R_n(\Gamma_g):= \{(x_1,x_2,...,x_{2g})\in (GL_n)^{2g}|\,\,r=1 \}.
\end{equation}
For example, when $\Gamma_g$ is the fundamental group of a compact Riemann surface, then the relation is 
$$
r= [x_1,x_2]....[x_i,x_{i+1}]...[x_{2g-1},x_{2g}]
$$ where $[x_i,x_{i+1}]$ denotes the commutator of $x_i$ and $x_{i+1}$.

We recall the following result on these varieties, by Rapinchuk, Benyash--Krivetz and Chernousov;
\begin{theorem}\label{irreducible}
The representation variety $ R_n(\Gamma_g)$ is an absolutely irreducible $\Q$--rational variety of dimension
\begin{eqnarray*}
\rm{dim }R_n(\Gamma_g)& = & (2g-1)n^2+1 \m{ if  } g>1 \\
 &=& n^2+n \m{ if  } g=1
\end{eqnarray*}
\end{theorem}
\begin{proof}
See \cite[Theorem 1]{RBC}.
\end{proof}

\subsection{Equivariant Chow-K\"unneth projectors for $R_n(\Gamma_g)$}
Since $GL_n$ is a homogeneous space it has only algebraic cohomology. Hence the same is true for the variety $(GL_n)^{2g}$. Hence, by Lemma \ref{simpleprojectors}, we get Chow--K\"unneth projectors for $(GL_n)^{2g}$.
Furthermore these projectors remain invariant under the conjugation action by
$Gl_n$ acting diagonally on the product and give equivariant Chow--K\"unneth projectors for $(GL_n)^{2g}$.

Let $R_n(\Gamma_g)^0$ be the smooth open subset of the variety $R_n(\Gamma_g)$.
This is an irreducible complex variety, by Theorem \ref{irreducible}.
Let $X_n(\Gamma_g)^0$ be the good quotient of
$R_n(\Gamma_g)^0$ (perhaps a smaller open subset so that quotient exists) modulo the action of $GL_n$ by conjugation.

\begin{theorem}
There is an equivariant Chow--K\"unneth decomposition for the smooth variety
$R_n(\Gamma_g)^0$, for $g>1$. This gives absolute Chow--K\"unneth projectors for the quotient variety  $X_n(\Gamma_g)^0$.
\end{theorem}
\begin{proof}
Using the description of the representation variety $R_n(\Gamma_g)$ in \eqref{relation}, the relation $r$ corresponds to $n^2-1$ equations in $(GL_n)^{2g}$ and hence $R_n(\Gamma_g)$ is a complete intersection in the product space
  $(GL_n)^{2g}$, defined by $n^2-1$ equations. 
There is a Lefschetz theorem for the smooth locus $R_n(\Gamma_g)^0$ (see \cite{D3}), and the restriction map
$$
H^i((GL_n)^{2g},\Q)\rar W_iH^i(R_n(\Gamma_g)^0,\Q)
$$
is bijective, for $i\leq \rm{dim }R_n(\Gamma_g)-1$.
Hence, by Lemma \ref{simpleprojectors}, we get equivariant Chow--K\"unneth 
projectors for the smooth variety $R_n(\Gamma_g)^0$. Using \cite{EdidinGraham} 
one can identify the rational equivariant Chow groups of $R_n(\Gamma_g)^0$ with the rational Chow groups of the quotient space $X_n(\Gamma_g)^0$, and this identification respects the ring structure. This means that the equivariant Chow--K\"unneth projectors
of $R_n(\Gamma_g)^0$ correspond to the absolute Chow--K\"unneth projectors for the smooth variety $X_n(\Gamma_g)^0$.
\end{proof}

Consider the commutator variety $C(2,GL_n)$ defined in \eqref{character}. 
Then, by \cite[Proposition 3]{RBC}, the variety $C(2,GL_n)$ is an absolutely 
irreducible $\Q$-rational variety.
Let $C(2,GL_n)^0\subset C(2,GL_n)$ be the smooth open locus. Then $C(2,GL_n)^0$ is a smooth irreducible complex variety.

\begin{theorem}
The smooth variety $C(2,GL_n)^0$ has an equivariant Chow--K\"unneth decomposition.
\end{theorem}
\begin{proof}
Firstly, using the description in \eqref{character}, we identify the variety $C(2g,GL_n)$ with the variety defined by quadratic relations in the product space $(GL_n)^{2g}$. In particular, when $g=1$, there is only one quadratic relation and $C(2,GL_n)^0$ is a smooth complete intersection in $(GL_n)^{2}$.
There is a Lefschetz theorem for the smooth locus $C(2,GL_n)^0$ (see \cite{D3}), and the restriction map
$$
H^i((GL_n)^{2},\Q)\rar W_iH^i(C(2,GL_n)^0,\Q)
$$
is bijective, for $i\leq \rm{dim }C(2,GL_n)-1$.

By applying Lemma \ref{simpleprojectors}, we get equivariant Chow--K\"unneth projectors for  $C(2,GL_n)^0$.  
\end{proof}


\end{document}